\documentclass[12pt,oneside]{amsart}
\usepackage{latexsym}
\usepackage{amsmath,amssymb,amsthm}
\usepackage{amscd}
\usepackage[dvips]{graphics}
\newcommand{\rp}{\mathbb{RP}}
\newcommand{\re}{\mathbb{R}}
\newcommand{\co}{\mathbb{C}}

\newcommand{\pgl}[1]{\mathbf{PGL}(#1,\mathbb{R})}
\newcommand{\Sl}[1]{\mathbf{SL}(#1,\mathbb{R})}
\newcommand{\gl}[1]{\mathbf{GL}(#1,\mathbb{R})}

\newcommand{\sfrac}[2]{{\textstyle \frac{#1}{#2}}}

\newtheorem{prop}{Proposition} 
\newtheorem{cor}[prop]{Corollary}
\newtheorem{thm}{Theorem}

\newtheorem{lem}[prop]{Lemma}

\theoremstyle{remark}

\begin{document}
\author{John Loftin}
\title{Flat Metrics, Cubic Differentials and Limits of Projective Holonomies}
\maketitle

\section{Introduction}
In \cite{labourie06} and \cite{loftin01}, it is shown that on a
closed oriented surface $S$ of genus $g>1$, there is a one-to-one
correspondence between   convex $\rp^2$ structures on $S$ and pairs
$(\Sigma,U)$, where $\Sigma$ is a conformal structure on $S$ and $U$
is a holomorphic cubic differential.  In this note, we compute the
asymptotic values of the holonomy of the $\rp^2$ structure
corresponding to $(\Sigma,\lambda U_0)$ as $\lambda\to\infty$ around
geodesic loops of the flat metric $|U_0|^\frac23$ which do not touch
any zeros of the fixed cubic differential $U_0$.  Such asymptotic
holonomies are related to the compactification of the deformation
space of convex $\rp^2$ structures on $S$ due to Inkang Kim
\cite{inkang-kim05} (see Section \ref{bound-def} below).
\begin{thm} \label{main-thm}
Let $\Sigma$ be a closed Riemann surface of genus $g>1$ and let
$U_0$ be a holomorphic cubic differential on $\Sigma$.  Consider a
closed oriented geodesic $\mathcal L$ of the flat metric
$|U_0|^\frac23$ on $\Sigma$ which does not touch any of the zeros of
$U_0$. In terms of the flat coordinate $z$ in which $U_0=2\,dz^3$,
represent the deck transformation corresponding to $\mathcal L$ as a
displacement $z\mapsto z+ Le^{i\theta}$ for $L>0$. Then there is a
constant $\kappa>0$ so that the eigenvalues $\xi_1>\xi_2>\xi_3$ of
the $\Sl3$ holonomy along $\mathcal L$ for the $\rp^2$ structure
determined by the pair $(\Sigma,\lambda U_0)$ for $\lambda>0$
satisfy
$$ \kappa \xi_i > e^{\lambda^\frac13\mu_i L} > \kappa^{-1} \xi_i$$
for $\mu_1\ge\mu_2\ge\mu_3$ the roots of the equation
$$\mu^3 - 3\mu - 2\cos3\theta = 0.$$
\end{thm}
  The techniques involved in the proof are similar
to the analysis of the harmonic map equation between hyperbolic
surfaces, as discussed by Mike Wolf \cite{wolf89} and Z.C.\ Han
\cite{han96}, and some new results on asymptotics of linear systems
of ODEs.

The present paper may be thought of as something of sequel to
\cite{loftin02c}, which studies the behavior of $\rp^2$ surfaces
corresponding to $(\Sigma,U)$ as $\Sigma$ approaches the boundary of
the Deligne-Mumford compactification of the moduli space of Riemann
surfaces, and $U$ degenerates to a regular cubic differential.

In future work, we hope to extend this analysis to all geodesics
with respect to the singular flat metric $|U_0|^\frac23$, including
those which are singular at the zeros of $U_0$.  This will allow a
full description of the data Kim prescribes for the boundary of the
deformation space of convex $\rp^2$ structures.  It will also be
interesting to relate the present work to harmonic maps to
$\re$-buildings, as an extension of Wolf's work on harmonic maps to
$\re$-trees \cite{wolf95}.

I would like to thank Mike Wolf, for, some years ago, pointing out
the similarities between the analytic theories of convex $\rp^2$
structures and harmonic maps between hyperbolic surfaces. I also
thank Bill Goldman for his encouragement and many fruitful
discussions about $\rp^2$ structures, and Lee Mosher for useful
discussions. The author is partially supported by NSF Grant
DMS0405873.

\section{The boundary of the deformation space of convex $\rp^2$
structures} \label{bound-def}

It is well known that a closed hyperbolic surface is determined by
its \emph{length spectrum}, which consists of the hyperbolic lengths
of the unique geodesic in each free homotopy class of curves. More
concretely, hyperbolic lengths of geodesic provide an embedding of
Teichm\"uller space into $\re^\mathcal C$, where $\mathcal C$ is the
set of all nontrivial conjugacy classes in $\pi_1(S)$ for a closed
surface $S$ of genus $g>1$.  Then Thurston's boundary of
Teichm\"uller space can be recovered as the set of limit points of
sequences in Teichm\"uller space $\subset \re^\mathcal C$, when
projected to the projective space $\mathbb {PR}^\mathcal C$
\cite{morgan-shalen84,bestvina88,paulin88}.

There is an analog of this theory to convex $\rp^2$ surfaces due to
Paulin \cite{paulin97}, Parreau \cite{parreau00} and Inkang Kim
\cite{inkang-kim01,inkang-kim05} (these authors address more general
structures as well). Recall a (properly) convex $\rp^2$ surface $S$
is given by $S=\Omega/\Gamma$, where $\Omega\Subset\re^2 \subset
\rp^2$ is a convex set and $\Gamma\subset\pgl3$.  Then each element
$\gamma\in\Gamma$ may be represented as a matrix in $\Sl3$. The
eigenvalues of this matrix are then analogs of the hyperbolic length
(see in particular Goldman \cite{goldman90a} for a detailed analog
of the Fenchel-Nielsen theory of Teichm\"uller space for the case of
convex $\rp^2$ structures).  In particular, for a given
$\gamma\in\Gamma$ with eigenvalues $\nu_1>\nu_2>\nu_3>0$ (the
eigenvalues have this structure by \cite{kac-vinberg67}), the set of
logarithms $$(\ell_1, \ell_2, \ell_3) = (\log \nu_1, \log \nu_2,
\log \nu_3)$$ is naturally an element of the maximal torus $\frak t$
of the Lie algebra $\frak{sl}(3,\mathbb R)$. Kim \cite{inkang-kim01}
shows that a normalized version of map of logarithms of eigenvalues
(into $\frak t^\mathcal C$) determines the $\rp^2$ structure. (The
normalization is an analog of projectivization of $\re^\mathcal C$
mentioned above for Teichm\"uller space.)  Then in
\cite{paulin97,parreau00,inkang-kim05}, the boundary of the
deformation space of convex $\rp^2$ structures may be defined to be
the boundary in  $\frak t^\mathcal C$ of the image of the
deformation space of all convex $\rp^2$ structures.

The limiting spectra in the case of both hyperbolic lengths and
$\rp^2$ structures can be seen as naturally arising in the context
of $\pi_1(S)$ actions on $\re$-buildings.  In particular, given a
background conformal structure $\Sigma_0$ on $S$, Teichm\"uller
space can be parametrized by the unique harmonic map from $\Sigma_0$
to the target hyperbolic structure \cite[Wolf]{wolf89}.  In turn,
these harmonic maps are uniquely determined by a holomorphic
quadratic differential $\Psi$ on $\Sigma_0$.  The key equation to
solve to construct the harmonic map is
$$ \Delta v + 4e^{-v}\|\Psi\|^2 - 2e^v + 2 = 0,$$ where $\Delta$ and
$\|\cdot\|$ are determined by the hyperbolic metric on $\Sigma_0$.
Wolf then essentially studies solutions to this equation to
reproduce Thurston's compactification of Teichm\"uller space as
limits of hyperbolic structures for quadratic differentials
$\lambda\Psi_0$ as $\lambda\to\infty$ \cite{wolf89}, and later uses
the same estimates to produce a $\pi_1(S)$ equivariant harmonic map
to an appropriate $\re$-tree \cite{wolf95}.

The equation we use (due to C.P.\ Wang \cite{wang91}) to produce
$\rp^2$ structures,
$$ \Delta u + 4e^{-2u} \|U\|^2 - 2e^u - 2\kappa = 0,$$ is
very similar to Wolf's equation, with a cubic differential $U$
replacing the quadratic differential $\Psi$.  There should also be
an analogous theory: The limiting structures for $U=\lambda U_0$ for
$\lambda\to\infty$ should be realized as an action of $\pi_1(S)$ on
a $\re$-building, together with a $\pi_1(S)$ equivariant map.  We
hope to address this problem in future work.

\section{Hyperbolic affine spheres and convex $\rp^n$ structures}
\label{h-a-s}

Recall the standard definition of $\rp^n$ as the set of lines
through 0 in $\re^{n+1}$. Consider
$\pi\!:\re^{n+1}\setminus0\to\rp^n$ with fiber $\re^*$.  For a
convex domain $\Omega\subset\re^n\subset\rp^n$ as above, then
$\pi^{-1}(\Omega)$ has two connected components. Call one such
component $\mathcal{C}(\Omega)$, the \emph{cone over $\Omega$}. Then
any representation of a group $\Gamma$ into $\pgl{n+1}$ so that
$\Gamma$ acts discretely and properly discontinously on $\Omega$
lifts to a representation into
$$\mathbf{SL}^\pm(n+1,\re)
= \{A\in\mathbf{GL}(n+1,\re):\det A=\pm1 \}$$ which acts on
$\mathcal{C}(\Omega)$. See e.g.\ \cite{loftin01}.

For a properly convex $\Omega$, there is a unique hypersurface
asymptotic to the boundary of the cone $\mathcal{C}(\Omega)$ called
the hyperbolic affine sphere
\cite{calabi72,cheng-yau77,cheng-yau86}. This hyperbolic affine
sphere $H\subset\mathcal{C}(\Omega)$ is invariant under
automorphisms of $\mathcal{C}(\Omega)$ in
$\mathbf{SL}^\pm(n+1,\re)$. The projection map $P$ induces a
diffeomorphism of $H$ onto $\Omega$. Affine differential geometry
provides $\mathbf{SL}^\pm(n+1,\re)$-invariant structure on $H$ which
then descends to $M=\Omega/\Gamma$.  In particular, both the affine
metric, which is a Riemannian metric conformal to the (Euclidean)
second fundamental form of $H$, and a projectively flat connection
whose geodesics are the $\rp^n$ geodesics on $M$, descend to $M$.
See \cite{loftin01} for details. A fundamental fact about hyperbolic
affine spheres is due to Cheng-Yau \cite{cheng-yau86} and
Calabi-Nirenberg (unpublished):

\begin{thm} \label{cheng-yau-thm}
If the affine metric on a hyperbolic affine sphere $H$ is complete,
then $H$ is properly embedded in $\re^{n+1}$ and is asymptotic to a
convex cone $\mathcal C\subset \re^{n+1}$ which contains no line.
By a volume-preserving affine change of coordinates in $\re^{n+1}$,
we may assume $\mathcal{C}=\mathcal{C}(\Omega)$ for some properly
convex domain $\Omega$ in $\rp^n$.
\end{thm}

Note that if $S=\Omega/\Gamma$ is compact, then Cheng-Yau's
completeness condition on the affine metric is satisfied by
\emph{any} appropriate affine metric on $S$.

\section{Wang's developing map}

In dimension 2, there is a local theory due to C.P.\ Wang
\cite{wang91} which exploits the elliptic PDE nature of the problem
of finding hyperbolic affine spheres to relate oriented convex
$\rp^2$ surfaces to holomorphic data on Riemann surfaces. See also
Labourie \cite{labourie06} and Loftin \cite{loftin02c}.  In
particular, the affine metric of a 2-dimensional hyperbolic affine
sphere induces a conformal structure on the surface, and, moreover,
there is a holomorphic cubic differential $U$ (which is essentially
the difference between the Levi-Civita connection of the affine
metric and the projectively flat connection of the $\rp^2$
structure) induced by the affine sphere.  All this structure
descends to projective quotients of the hyperbolic affine sphere. In
particular, we have the following
\begin{thm}
Given an oriented surface $S$, the structure of a convex $\rp^2$
structure on $S$ is equivalent to the pair of a conformal structure
$\Sigma$ and a holomorphic cubic differential $U$ on $S$.
\end{thm}

Locally, the structure equations of a 2-dimensional hyperbolic
affine sphere may be expressed in terms of a embedding map
$f\!:\Omega \to \re^3$, where $\Omega\subset \co$ is a
simply-connected domain.  $f$ is taken to be a conformal map with
respect to the affine metric $e^\psi|dz|^2$ and $U$ is a holomorphic
function.  Then $f$ satisfies

\begin{equation}
\left\{ \begin{array}{c}
f_{zz} = \psi_z f_z + U e^{-\psi} f_{\bar z} \\
f_{{\bar z}{\bar z}} = {\bar U} e^{-\psi} f_z + \psi_{\bar z}
f_{\bar z}
\\
f_{z{\bar z}} = \frac{1}{2}e^\psi f \end{array} \right.
\label{fzz-eq}
\end{equation}
The conformal factor $e^\psi$ must satisfy the following
integrability condition,
\begin{equation}
\psi_{z {\bar z}} + |U|^2 e^{-2\psi} - \sfrac{1}{2} e^\psi  =0,
\label{loc-wang-eq}
\end{equation}
which we call Wang's equation.  On a Riemann surface $U$ transforms
as a cubic differential, and (\ref{loc-wang-eq}) becomes, with
respect to a conformal background metric $h$,
\begin{equation}
\Delta u + 4 e^{-2u} \|U\|^2 - 2 e^u - 2\kappa=0, \label{wang-eq}
\end{equation}
where $\Delta$ is the Laplacian of $h$, $\|U\|^2$ is the
norm-squared of $U$ with respect to the metric $h$, $\kappa$ is the
Gauss curvature of $h$, and the metric $e^uh = e^\psi|dz|^2$ locally
for $\psi$ given by (\ref{loc-wang-eq}).

We now study solutions to (\ref{wang-eq}) for $U=\lambda U_0$ as
$\lambda\to\infty$.

\section{Limits of the conformal metrics}
Let $U_0$ be a holomorphic cubic differential on $\Sigma$ which is
not identically zero.  We study the limiting behavior of solutions
to Wang's equation (\ref{wang-eq}) for solutions $u_{\lambda}$ as $U
= \lambda U_0$ for $\lambda$ a real parameter approaching $\infty$.
In his work on harmonic maps between hyperbolic surfaces and
Thurston's boundary of Teichm\"uller space, Mike Wolf has studied a
similar equation to (\ref{wang-eq}) with a holomorphic quadradic
differential instead of a cubic differential \cite{wolf89}. The
proof below is similar to the one in Han \cite{han96}.

\begin{prop} \label{u-asymp}
Let $\Sigma$ be a closed Riemann surface of genus $g>1$ equipped
with a background metric $h$ and a holomorphic cubic differential
$U_0$ which is not identically zero.  Let $\lambda>0$ and let
$u=u_\lambda$ be the solution to (\ref{wang-eq}) for $U=\lambda
U_0$. Let $K$ be a compact subset of $\Sigma$ which does not contain
any of the zeroes of $U_0$. Then there is a constant
$C=C(\Sigma,U_0,K)$ so that
$$ \frac12 \ge \|U\|^2e^{-3u_\lambda} \ge \frac12 - C\lambda^{-\frac23}.$$
\end{prop}

\begin{proof}
We prove this proposition by the use of barriers.  The key
observation is that the singular flat conformal metric
$2^\frac13|U|^\frac23$ provides a solution to (\ref{loc-wang-eq})
away from the zeros of $U$.

Consider a smooth background metric $g$ by requiring
$g=2^\frac13|U_0|^\frac23$ on $K$ and $\|U_0\|^2_g \le \frac12$ on
all $\Sigma$. (This is possible since $\|U_0\|^2_g=\frac12$ on $K$
and $\|U_0\|^2_g = 0$ at the zeros of $U_0$.)

Now for $U=\lambda U_0$, define $s=s_\lambda$ by
$$ge^s = 2^\frac13|U|^{\frac23}=2^\frac13\lambda^\frac23 |U_0|^\frac23.$$
Note that $s=\frac23\log\lambda$ on $K$. We may also check that $s$
solves (\ref{wang-eq}) away from the zeros of $U$, and is equal to
$-\infty$ at the zeros of $U$.  By applying the comparison principle
to (\ref{wang-eq}), we find that $u\ge s$, and so $s$ is a
subsolution of (\ref{wang-eq}).

Now let $S=S_\lambda$ be equal to $\log r$ for $r=r_\lambda$ the
positive root of $$p(x)=x^3 - \sigma x^2 -\lambda^2 = 0, \qquad
\sigma = \max_\Sigma (-\kappa_g),$$ for $\kappa_g$ the Gauss
curvature of $g$.  Then $S$ is a supersolution of (\ref{wang-eq}):
At a maximum point of $u$,
\begin{eqnarray*}0&\ge&\Delta_g u = 2e^u + 2\kappa -
4e^{-2u}\|U\|_g^2, \\
&\ge& 2e^{-2u}(e^{3u} - \sigma e^{2u} - \lambda^2).
\end{eqnarray*}
The largest value of $u$ for which this inequality can be true
occurs when $p(e^u)=0$.

On $K$ then, $\frac23\log\lambda\le u\le S$, and so
$$\frac12\ge\|U\|^2_g e^{-3u}\ge \frac12 \,\lambda^2e^{-3S}.$$
Now we note that $\tilde x =\lambda^{-\frac23}e^S$ solves
$$ \tilde x^3 - \sigma\lambda^{-\frac23}\tilde x^2 - 1 = 0,$$ and
so for large values of $\lambda$, $\tilde x = 1 +
O(\lambda^{-\frac23})$.  This proves the proposition.
\end{proof}

\begin{cor} \label{d-psi-estimate}
There is another constant $C=C(\Sigma,U_0,K)$ so that $|\psi_z|\le
C\lambda^{-\frac13}$ on $K$, where $z$ is a local coordinate so that
$U_0=2\,dz^3$.
\end{cor}

\begin{proof}
Note that in the proof and below, different uniform constants may be
referred to by the same letter $C$ depending on the context.

For $p\in K$, choose the local coordinate $z$ so that $z(p)=0$ and
let consider $$\alpha(w) = \psi(\lambda^{-\frac13}w)
-\sfrac23\log\lambda -\sfrac13\log 2.$$
 Then $$\alpha_{w\bar w} = \lambda^{-\frac23}\psi_{z\bar z} =
 2^{-\frac23}(e^{-2\alpha}-e^\alpha).$$
Proposition (\ref{u-asymp}) implies that there is a constant $C$ so
that
 $$0\le \alpha(\lambda^\frac13z) =\psi(z)-\sfrac23\log\lambda -\sfrac13\log 2
 \le C \lambda^{-\frac23}$$
for all $z$ in a neighborhood of $K$.

This implies that in any disk in the $w$-plane centered at 0, there
is a constant $C$ independent of $p\in K$ and $\lambda$ large so
that
$$\left|\alpha\right|, \left|\alpha_{w\bar w}\right| \le C
\lambda^{-\frac23}.$$  Then the $L^p$ theory implies that on a
slightly smaller disk, that $\|\alpha\|_{W^{2,p}} \le C
\lambda^{-\frac23}$.  Then, for $p>2$, Sobolev embedding implies
similar bounds for the $C^1$ norm of $\alpha$:
$$ |\alpha_w| \le C \lambda^{-\frac23}.$$
Now simply compute $\psi_z =\lambda^{\frac13} \alpha_w$.
\end{proof}

\section{ODE estimates}
Now the structure equations (\ref{fzz-eq}) can be recast in terms of
the frame $\langle f,\lambda^{-\frac13} f_z, \lambda^{-\frac13}
f_{\bar z}\rangle$ to read
\begin{eqnarray} \label{fzz}
\left(\begin{array}{c} f\\ \lambda^{-\frac13}f_z \\
\lambda^{-\frac13}f_{\bar z}\end{array}\right)_z &=& \left(
\begin{array}{ccc} 0&\lambda^{\frac13}&0 \\ 0& \psi_z & Ue^{-\psi}\\
\frac12\lambda^{-\frac13}e^\psi & 0& 0\end{array} \right)
\left(\begin{array}{c} f\\ \lambda^{-\frac13}f_z \\
\lambda^{-\frac13}f_{\bar z}\end{array}\right), \\
 \label{fzzbar}
\left(\begin{array}{c} f\\ \lambda^{-\frac13}f_z \\
\lambda^{-\frac13}f_{\bar z}\end{array}\right)_{\bar z} &=& \left(
\begin{array}{ccc} 0&0&\lambda^{\frac13} \\
\frac12\lambda^{-\frac13}e^\psi & 0& 0
\\ 0& \bar U e^{-\psi}& \psi_{\bar z}
\end{array} \right)
\left(\begin{array}{c} f\\ \lambda^{-\frac13}f_z \\
\lambda^{-\frac13}f_{\bar z}\end{array}\right).
\end{eqnarray}

Away from the zeros of $U_0$, choose a local coordinate $z$ so that
$U_0 = 2\,dz^3$, and $U=\lambda U_0 = 2\lambda\,dz^3$. Proposition
\ref{u-asymp} and Corollary \ref{d-psi-estimate} then show that the
matrices in the structure equations above have the form
\begin{eqnarray}
P&=&\left(
\begin{array}{ccc} 0&\lambda^{\frac13}&0 \\ 0& \psi_z & Ue^{-\psi}\\
\frac12\lambda^{-\frac13}e^\psi & 0& 0\end{array} \right) =
\lambda^\frac13\left(\begin{array}{ccc} 0&1&0 \\
0&0&1 \\
1&0&0 \\
\end{array}\right) + O(\lambda^{-\frac13}),\\
Q&=&\left(
\begin{array}{ccc} 0&0&\lambda^{\frac13} \\
\frac12\lambda^{-\frac13}e^\psi & 0& 0
\\ 0& \bar U e^{-\psi}& \psi_{\bar z}
\end{array} \right) =
\lambda^\frac13\left(\begin{array}{ccc} 0&0&1 \\
1&0&0 \\
0&1&0 \\
\end{array}\right) + O(\lambda^{-\frac13}),
\end{eqnarray}
where $O(\lambda^{-\frac13})$ is as $\lambda\to\infty$ for all
points in $K$ a compact set not containing any zero of $U_0$.

We will integrate the initial value problem along a geodesic path
with respect to the metric $|U_0|^\frac23$ which avoids the zeroes
of $U_0$.  These paths are simply straight lines each local complex
coordinate chart with coordinate $z$ satisfying $U_0=2\,dz^3$.  In
the particular case of a geodesic loop, the system of ODEs
(\ref{fzz}-\ref{fzzbar}) will compute the real projective holonomy
around such a loop: along such a loop, the coordinate $z$ can be
analytically continued in the universal cover, and the corresponding
deck transformation corresponds to $z\mapsto z+c$ for a complex
constant $c$.  Therefore, the frame $\langle f,\lambda^{-\frac13}
f_z, \lambda^{-\frac13} f_{\bar z}\rangle$ is the frame of a rank-3
vector bundle on the quotient whose holonomy in $\gl3$ projects to
$\pgl3$ to compute the real projective holonomy of the geodesic
loop.  For more details of this argument, see e.g.\ Proposition 2 of
\cite{loftin02c}.

Any geodesic loop of $|U_0|^\frac23$ which avoids the zeroes of
$U_0$ may be described by a starting point, at which we set the
local coordinate $z$ to be $0$, and a finishing point, which we set
to be $z=c$ in the analytically continued $z$ coordinate.  The
geodesic is then the straight line segment between 0 and $c$. If
$c=Le^{i\theta}$ for $L>0$, then the holonomy with respect to the
frame $\mathcal F = \langle f,\lambda^{-\frac13} f_z,
\lambda^{-\frac13} f_{\bar z}\rangle$ is $\Phi(L)$, where $\Phi$
solves the initial value problem $$ \Phi(0)=I, \qquad
\frac{d\Phi}{dt} =(e^{i\theta}P + e^{-i\theta}Q)\Phi.$$ This ODE
system is equivalent to
$$ \frac{d\Phi}{dt} = \left[\lambda^\frac13 \left(\begin{array}{ccc} 0&
e^{i\theta} & e^{-i\theta} \\ e^{-i\theta} &0&e^{i\theta} \\
e^{i\theta} & e^{-i\theta} &0 \end{array}\right) +
O(\lambda^{-\frac13}) \right] \Phi.$$ As we are primarily interested
in the eigenvalues of $\Phi(L)$, we replace the matrix
$$M=\left(\begin{array}{ccc} 0&
e^{i\theta} & e^{-i\theta} \\ e^{-i\theta} &0&e^{i\theta} \\
e^{i\theta} & e^{-i\theta} &0 \end{array}\right) $$ by the conjugate
diagonal matrix
$$\left(\begin{array}{ccc}\mu_1&0&0\\ 0&\mu_2&0 \\ 0&0&\mu_3
\end{array}\right)$$ for $\mu_i$ the roots of the characteristic
equation $$ \det(\mu I- M)=\mu^3-3\mu-2\cos3\theta=0.$$ We note $M$
is diagonalizable and $\mu_i\in\re$. Assume $\mu_1\ge\mu_2\ge\mu_3$.

Then, to compute the conjugacy class of the holonomy matrix around
this geodesic loop, we compute the solution to
\begin{eqnarray}
\label{phi-init-cond}
\Phi(0)&=&I,\\
\label{phi-ode}
 \frac{d\Phi}{dt} &=&
\left[\lambda^{\frac13}\left(\begin{array}{ccc} \mu_1&0&0 \\
0&\mu_2&0 \\ 0&0&\mu_3 \end{array}\right) +
\left(\begin{array}{ccc}b_{11}&b_{12}&b_{13} \\ b_{21} & b_{22} &
b_{23} \\ b_{31} & b_{32} & b_{33} \end{array}\right) \right] \Phi,
\end{eqnarray}
where there is a constant $C$ so that each $b_{ij} =
b_{ij}(t,\lambda)$ satisfies $|b_{ij}|\le C\lambda^{-\frac13}$.

\begin{prop} \label{fund-sol-asymp}
The solution $\Phi$ to the initial value problem
(\ref{phi-init-cond}-\ref{phi-ode}) has the form
$$\left(\begin{array}{ccc} e^{\lambda^\frac13\mu_1 t}+O(\lambda^{-\frac13}
e^{\lambda^\frac13\mu_1 t}) & O(\lambda^{-\frac13}
e^{\lambda^\frac13\mu_2 t})& O(\lambda^{-\frac13}
e^{\lambda^\frac13\mu_3 t})\\
O(\lambda^{-\frac13} e^{\lambda^\frac13\mu_1 t})&
e^{\lambda^\frac13\mu_2 t}+O(\lambda^{-\frac13}
e^{\lambda^\frac13\mu_2 t})& O(\lambda^{-\frac13}
e^{\lambda^\frac13\mu_3 t})\\
O(\lambda^{-\frac13} e^{\lambda^\frac13\mu_1 t})&
O(\lambda^{-\frac13} e^{\lambda^\frac13\mu_2 t})
&e^{\lambda^\frac13\mu_3 t}+O(\lambda^{-\frac13}
e^{\lambda^\frac13\mu_3 t})
\end{array}\right),$$
Where the $O$ notation denotes bounds as $\lambda\to\infty$ that are
uniform for $t\in[0,L]$.
\end{prop}

\begin{proof}
Write $\Phi = (\phi_{ij})$, and consider the first column
$\phi_{11}, \phi_{21},\phi_{31}$, which satisfies the linear system
\begin{eqnarray*}
 \phi_{11}(0)= 1, &\qquad&\sfrac d{dt}\phi_{11} =
 (\lambda^\frac13\mu_1+b_{11})\phi_{11} + b_{12}\phi_{21}+
b_{13}\phi_{31}, \\
 \phi_{21}(0)= 0, &\qquad &
\sfrac d{dt} \phi_{21} =
 b_{21}\phi_{11}
+(\lambda^{\frac13}\mu_2+b_{22})\phi_{21} + b_{23}\phi_{31}, \\
 \phi_{31}(0)=0,&\qquad&
\sfrac d{dt}\phi_{31} = b_{31}\phi_{11} + b_{32}\phi_{21} +
(\lambda^\frac13\mu_3 + b_{33}) \phi_{31}.
\end{eqnarray*}

Each of the above differential equations is first-order linear, and
so we must have
\begin{eqnarray*}
 \phi_{11} &=& e^{\lambda^\frac13\mu_1t}e^{\int_0^tb_{11}}
 \left[ 1+\int_0^t e^{-\lambda^\frac13 \mu_1 \tau-\int_0^\tau
 b_{11}} (b_{12}\phi_{21}+b_{13}\phi_{31})d\tau
 \right], \\
\phi_{21} &=& e^{\lambda^\frac13\mu_2t}e^{\int_0^tb_{22}}
 \int_0^t e^{-\lambda^\frac13 \mu_2 \tau-\int_0^\tau
 b_{22}} (b_{21}\phi_{11}+b_{23}\phi_{31})d\tau, \\
 \phi_{31} &=& e^{\lambda^\frac13\mu_3t}e^{\int_0^tb_{33}}
 \int_0^t e^{-\lambda^\frac13 \mu_3 \tau-\int_0^\tau
 b_{33}} (b_{31}\phi_{11}+b_{32}\phi_{21})d\tau.
\end{eqnarray*}
The previous three equations can be seen as a map $\mathcal M$ from
the $\re^3$-valued function $(\phi_{11},\phi_{21},\phi_{31})$ to the
right-hand sides.

Now let $N\gg1$ be a constant independent of $\lambda$, and consider
the Banach space $\mathcal B_\lambda$ of continuous $\re^3$-valued
functions with norm
$$\|(f_1,f_2,f_3)\|_{\mathcal B_\lambda}
 = \sup_i \sup_{t\in[0,L]} |f_i(t)| e^{-\lambda^\frac13
 \mu_1 t}.$$
Let $\mathcal B_\lambda(N)$ be the closed ball of radius $N$
centered at the origin in $\mathcal B_\lambda$. We claim that for
$\lambda$ large enough, $\mathcal M$ is a contraction map from
$\mathcal B_\lambda(N)$ to itself, and thus the solution
$(\phi_{11},\phi_{21},\phi_{31})$ to the ODE system, which is the
fixed point of $\mathcal M$, must lie in $\mathcal B_\lambda(N)$.

Now consider $F=(f_1,f_2,f_3),$ $G=(g_1,g_2,g_3)\in \mathcal
B_\lambda(N)$. Then the first component of $\mathcal M(F)-\mathcal
M(G)$ is given by
$$ e^{\lambda^\frac13 \mu_1 t} e^{\int_0^tb_{11}} \int_0^t
e^{-\lambda^\frac13 \mu_1 \tau - \int_0^\tau b_{11}}[b_{12}
(f_2-g_2) + b_{13}(f_3-g_3)]d\tau.$$ Now assume $|b_{ij}|\le R$ and
recall $t\le L$. Then a straightforward calculation shows that the
first component of $\mathcal M(F)-\mathcal M(G)$ is pointwise
bounded by
$$e^{\lambda^\frac13 \mu_1 t} e^{2RL}\cdot 2R \cdot L\cdot \|F-G\|_
{\mathcal B_{\lambda}},$$ and so if we choose $\lambda$ large enough
so that $R\sim \lambda^{-\frac13}$ is small enough, we may assume
$e^{2RL}\cdot 2R\cdot L <1$. Essentially the same calculation shows
that $\mathcal M\!:\mathcal B_{\lambda}(N) \to \mathcal
B_\lambda(N)$ for large $\lambda$, since $N\gg1$.  The two other
components of $\mathcal M$ behave the same way. All this shows
$\mathcal M$ is a contraction map.

Since $\mathcal M$ is a contraction map on the complete metric space
$\mathcal B_\lambda(N)$, the unique solution
$(\phi_{11},\phi_{21},\phi_{31})$ to the ODE system is the fixed
point, and so must be in $B_\lambda(N)$ for all $\lambda$
sufficiently large. Now simply apply the bounds
$$|\phi_{11}|, \, |\phi_{21}|, \, |\phi_{31}| \le N
e^{\lambda^\frac13 \mu_1 t}$$ to the fixed point equation
$(\phi_{11},\phi_{21},\phi_{31}) = \mathcal M(\phi_{11}, \phi_{21},
\phi_{31})$ to show that
$$\phi_{11} = e^{\lambda^\frac13 \mu_1 t} +
O(\lambda^{-\frac13}e^{\lambda^\frac13 \mu_1 t}), \quad \phi_{21} =
O(\lambda^{-\frac13}e^{\lambda^\frac13 \mu_1 t}),\quad \phi_{31} =
O(\lambda^{-\frac13}e^{\lambda^\frac13 \mu_1 t}).$$

 This justifies the first column in the matrix in Proposition
 \ref{fund-sol-asymp}. The argument for the other two columns is
 identical.
\end{proof}

\begin{thm}
There is a constant $\kappa>0$ so that the eigenvalues $\xi_1\ge
\xi_2\ge \xi_3>0$ of the holonomy matrix $\Phi(L)$ satisfy
$$\kappa\xi_i
> e^{\lambda^\frac13\mu_iL} > \kappa^{-1}\xi_i$$
for $i=1,2,3$.
\end{thm}

\begin{proof}
Proposition \ref{fund-sol-asymp} and the fact that $\Phi(L)\in \Sl3$
show that the characteristic polynomial of $\Phi(L)$ is $$
\begin{array}{c} x^3 - (e^{\lambda^\frac13\mu_1L} +
e^{\lambda^\frac13\mu_2L} +
e^{\lambda^\frac13\mu_3L})[1+O(\lambda^{-\frac13})]x^2 +{} \\
(e^{\lambda^\frac13(\mu_1+\mu_2)L} + e^{\lambda^\frac13(\mu_1
+\mu_3)L} + e^{\lambda^\frac13(\mu_2+\mu_3)L})
[1+O(\lambda^{-\frac13})]x - 1. \end{array}$$ Kac-Vinberg showed
\cite{kac-vinberg67} that the holonomy of any nontrivial loop in a
closed oriented convex $\rp^2$ surface of genus $g>1$ has positive
distinct eigenvalues $\xi_1>\xi_2>\xi_3>0$. Then
 $$\xi_1+\xi_2 + \xi_3 = (e^{\lambda^\frac13\mu_1L} +
e^{\lambda^\frac13\mu_2L} +
e^{\lambda^\frac13\mu_3L})[1+O(\lambda^{-\frac13})]$$ implies that
there is an $\epsilon$ which goes to 0 as $\lambda\to\infty$ so that
$$ (3+\epsilon)e^{\lambda^\frac13\mu_1L} > \xi_1 > (\sfrac13-\epsilon)
e^{\lambda^\frac13\mu_1L}.$$ Now use the bounds on $\xi_1$ and
$$\xi_1\xi_2+\xi_1\xi_3+\xi_2\xi_3 =
(e^{\lambda^\frac13(\mu_1+\mu_2)L} + e^{\lambda^\frac13(\mu_1
+\mu_3)L} + e^{\lambda^\frac13(\mu_2+\mu_3)L})
[1+O(\lambda^{-\frac13})]$$ to conclude that there is an
$\epsilon'\to0$ as $\lambda\to\infty$ so that
$$ (9+\epsilon')e^{\lambda^\frac13\mu_2L}  > \xi_2 >
(\sfrac19-\epsilon')e^{\lambda^\frac13\mu_2L}.$$ Then the theorem
follows from $\mu_1+\mu_2+\mu_3=0$ and  $$\xi_1\xi_2\xi_3=1.$$
\end{proof}

Since $\Phi(L)$ is conjugate to the holonomy matrix with respect to
the frame $\langle f, \lambda^{-\frac13}f_z, \lambda^{-\frac13}
f_{\bar z} \rangle$ around the loop $\mathcal L$, this concludes the
proof of Theorem \ref{main-thm}.

\bibliographystyle{abbrv}
\bibliography{thesis}

\end{document}